\documentclass[12pt]{article}
\usepackage{graphicx}
\usepackage{latexsym,array,delarray,amsthm, amssymb}
\usepackage[sumlimits]{amsmath}

\theoremstyle{plain}
\newtheorem{thm}{Theorem}[section]
\newtheorem{lemma}[thm]{Lemma}
\newtheorem{prop}[thm]{Proposition}
\newtheorem{cor}[thm]{Corollary}
\newtheorem{conj}[thm]{Conjecture}
\theoremstyle{definition}

\newtheorem{ex}[thm]{Example}

\newtheorem{algorithm}[thm]{Algorithm}
\theoremstyle{remark}

\newcommand{\rr}{\mathbb{R}}

\newcommand{\cB}{\mathcal{B}}
\newcommand{\cF}{\mathcal{F}}

\begin{document}

\title{Least Squares Methods for Equidistant Tree Reconstruction}
\author{Conor Fahey, Serkan Ho\c{s}ten, \\
Nathan Krieger, Leslie Timpe \\
{\small Department of Mathematics, San Francisco State University} }
\date{}
\maketitle

\begin{abstract}
UPGMA is a heuristic method identifying the least squares equidistant
phylogenetic tree given empirical distance data among $n$ taxa. We study
this classic algorithm using the geometry of the space of all equidistant
trees with $n$ leaves, also known as the Bergman complex of the graphical
matroid for the complete graph $K_n$. We show that UPGMA performs an
orthogonal projection of the data onto a maximal cell of the Bergman
complex. We also show that the equidistant tree with the least (Euclidean)
distance from the data is obtained from such an orthogonal projection, but
not necessarily given by UPGMA. Using this geometric information we give an
extension of the UPGMA algorithm. We also present a branch and bound method
for finding the best equidistant tree. Finally, we prove that there are
distance data among $n$ taxa which project to at least $(n-1)!$ equidistant
trees.

\end{abstract}


\section{Introduction}

\label{sec:intro}

\noindent We study the problem of finding the least squares equidistant tree
given distance data between the elements of a finite set $X$ of cardinality $%
n$. The set $X$ is often a collection of taxa in biological applications. In
this paper we will usually let $X=\{1,2,\ldots ,n\}$ unless otherwise
stated. The distance data is given by a \emph{dissimilarity map}, a
real-valued function $d:{\binom{[n]}{2}}\longrightarrow \mathbb{R}$ defined
for pairs $(i,j)$ where $1\leq i<j\leq n$. We will represent a dissimilarity
map by the edge weights of the complete graph $K_{n}$ of $n$ vertices where
the weight of the edge $(i,j)$ in $K_{n}$ is $d(i,j)$.

\noindent Let $T$ be a (not necessarily binary) weighted tree with a root $r$
and $n$ leaves. The weight of each edge $e\in T$ will be denoted by $%
w_{T}(e) $, and we will omit the subscript when the context makes it clear
which tree we refer to. Given such a tree we get a distance function $%
x(a,b)=\sum_{e\in P_{a,b}}w_{T}(e)$ where $P_{a,b}$ is the unique path
between the nodes $a$ and $b$ in $T$. A tree is called \emph{equidistant} if 
$x(i,r)$ is the same real number for each leaf $i=1,\ldots ,n$. Note that we
label the leaves of $T$ by $X$. A tree is equidistant if and only if for
each distinct $i,j,k\in X,$ the set of distances $\{x(i,j),x(i,k),x(j,k)\}$
achieves its maximum at least twice \cite{Bu,Sim,Z}. These are the \textit{%
ultrametric} conditions.

\vskip0.2cm \noindent With these definitions we can present the main problem
of this paper: given a dissimilarity map $d$ on $X=\{1,\ldots ,n\}$ find an
equidistant tree $T$ on $n$ leaves such that 
\begin{equation*}
\sum_{1\leq i<j\leq n}(d(i,j)-x(i,j))^{2}
\end{equation*}%
is minimized. It is known that this problem (as well as the unrooted
nonequidistant version) is NP complete \cite{K,KM84,KM86}.

\vskip0.2cm \noindent The problem of tree construction arises in biology,
where the goal is to describe the evolutionary history of species or genes.
An equidistant tree approximates the true evolutionary history. The
distances between species may be measured using several different methods,
but currently distances are most often determined by comparison of aligned
nucleic acid or amino acid sequences. One of several models of evolution is
used to correct for the possibility of multiple substitutions at any one
site \cite{Fel}. When the rate of nucleotide or amino acid substitution was
constant over the time period being considered, the ultrametric conditions
are close to being satisfied. This condition is the \emph{molecular clock
hypothesis,} and if it holds a least squares equidistant tree could be used
to fit the distance data. Least squares methods for tree construction are
attractive because they are statistically consistent: \ the correct tree
will be identified in the limit as the length of the sequences grows \cite%
{Fel84,Fel}. In many cases the molecular clock hypothesis is not
satisfied, and trees that are additive but not equidistant are preferred. \ 

\noindent
The Unweighted Pair Group Method with Arithmetic Means (UPGMA) algorithm is
a heuristic method for finding the least squares equidistant tree \cite{Fel}. \
The UPGMA algorithm has polynomial time complexity, and works well on data
which shows clock-like behavior. \ Even if the molecular clock holds,
however, the UPGMA algorithm may return a tree that is not the best by the
least squares criterion, as shown in Example 2.5 below. \ The unweighted
least squares approach was first suggested by Cavalli-Sforza and Edwards \cite{CE}. 
\ Other, related algorithms include the pioneering weighted least
squares algorithm of Fitch and Margoliash \cite{FM}, the transformed
distances method \cite{Fa}, and neighbor-joining \cite{SN}.
Neighbor-joining and recent variants BIONJ \cite{Ga} \ and weighbor \cite{BSH}
are not strictly least squares algorithms. \ Of all these UPGMA is
particularly interesting here, as it arises naturally as a greedy algorithm
from the approach described below. When the Euclidean metric is replaced by $%
\ell _{\infty }$ metric, a fast exact algorithm is known \cite{CF}. A
conceptual explanation of this algorithm is given in \cite{A}. \ 

Here we first describe UPGMA. We will present a version that outputs the
combinatorial description of $T$ and $x(i,j)$ for each pair of leaves $i$
and $j$. It is well-known how to compute the edge weights $w_{T}(e)$ from
these data. Recall that we represent the dissimilarity map $d$ as the edge
weights of $K_{n}$.

\begin{algorithm}
\label{alg:upgma} UPGMA \newline
\textbf{Input} : Complete graph $K_n$ with edge weights $d(i,j)$. \newline
\textbf{Output}: An equidistant tree $T$ with leaves $X = \{1, \ldots, n\}$
and $x(i,j)$ for each $i,j \in X$. \newline

\vskip 0.3cm \noindent $G := K_n$. \newline
$V(T) := X$, $E(T) := \emptyset$, and $S(T) := X$. \newline
\textbf{repeat} \newline
\begin{equation*}
minave := \min_{v,w \in V(G)} \frac{1}{C(v,w)} \sum_{(i,j) \in E(v,w)} d(i,j)
\end{equation*}
where $E(v,w)$ is the set of edges between the nodes $v$ and $w$ in $G$, and 
$C(v,w) = |E(v,w)|$. \newline
Let $s$ and $t$ in $V(G)$ be the vertices for which the minimum above is
attained. \newline
Set $x(i,j) := minave$ for all $(i,j) \in E(s,t)$. \newline
$G := G / \{s,t\}$, obtained by contracting the vertices $s$ and $t$ into a
single vertex $st := s \cup t$. \newline
$S(T) := S(T) \setminus \{s, t\} \cup \{st \}$. \newline
$V(T) := V(T) \cup \{st\}$, $E(T) := E(T) \cup \{st, s\} \cup \{st, t\}$. 
\newline
\textbf{until} $G$ has one vertex \newline

\vskip 0.2cm \noindent Output $T$ and $x(i,j)$ $1 \leq i < j \leq n$. \hfill 
$\Box$
\end{algorithm}

\noindent In Section \ref{sec:bergman} we describe the Bergman complex $%
\mathcal{B}_{n}$, namely, the space of all equidistant trees on $n$ leaves.
We prove that Algorithm \ref{alg:upgma} performs an orthogonal projection
(with respect to the usual Euclidean inner product) onto a maximal cell of $%
\mathcal{B}_{n}$. We give an example where already for $n=4$ the UPGMA tree
can be arbitrarily worse than the best equidistant tree. In Section \ref%
{sec:projections} we prove that the best equidistant tree is obtained by an
orthogonal projection onto some maximal cone of $\mathcal{B}_{n}$. Motivated
by this result we introduce a polyhedral subdivision of the data space $%
\mathbb{R}^{\binom{n}{2}}$ where each maximal cell consists of data vectors
which project onto the same set of maximal Bergman cells. We show that the
collection of such Bergman cells could be disconnected (in a sense made
precise in Section \ref{sec:projections}). In fact, there are data vectors in $\mathbb{%
R}^{\binom{n}{2}}$ which project onto at least $(n-1)!$ Bergman cells, and
we conjecture that this is the most number of projections one can obtain.
Furthermore, we classify all data vectors in $\mathbb{R}^{6}$ which project
onto six Bergman cells. In Section \ref{sec:comp} we introduce two
algorithms based on our results in Section \ref{sec:projections}. One of
them is an extension of UPGMA that finds at least as good a tree as the
UPGMA tree. The other one finds the best equidistant tree using a branch and
bound approach. Section \ref{sec:biology} concludes with an example where we
analyze data for the timing and the sequence of the appearance of mammalian
orders.


\section{The Bergman complex and UPGMA}

\label{sec:bergman}

\noindent It is not difficult to show that Algorithm \ref{alg:upgma} indeed
returns an equidistant tree using the ultrametric characterization of
equidistant trees. Here we will describe the space of all vectors $%
x=(x(i,j)\,:\,1\leq i<j\leq n)\in \mathbb{R}^{\binom{n}{2}}$ which come from
weighted equidistant trees with $n$ leaves, and from this description it
will follow that the UPGMA produces an equidistant tree. Ardila and Klivans 
\cite{AK} described this space as a special case of the tropicalization of a
linear variety, or more combinatorially, as the \emph{Bergman complex} $%
\mathcal{B}_{n}$ of the graphical matroid of $K_{n}$. This description shows
that $\mathcal{B}_{n}\subset \mathbb{R}^{\binom{n}{2}}$ is a polyhedral
complex of dimension $n-1$: its maximal cones are polyhedral cones of
dimension $n-1$, and any collection of them intersects in a face that
belongs to each cone in the collection.

\begin{figure}[tbp]
\centerline{
\includegraphics[height=10cm]{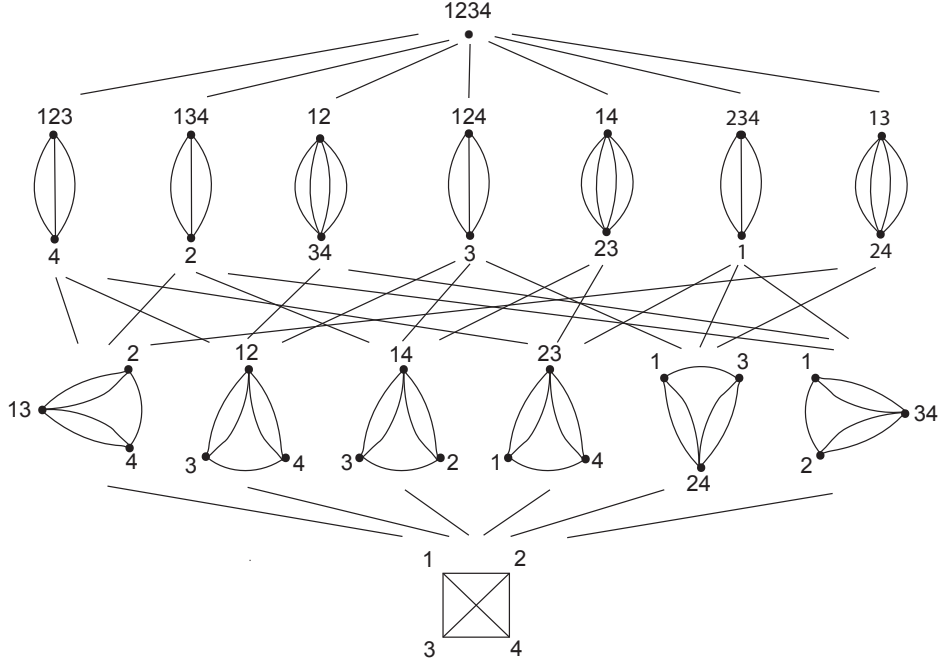}
}
\caption{Lattice of contractions of $K_4$}
\label{Figure1}
\end{figure}

\noindent We first describe a different polyhedral complex $\mathcal{F}_{n}$
of dimension $n-1$ that is a refinement of $\mathcal{B}_{n}$, i.e. the
maximal cones of $\mathcal{F}_{n}$ further subdivide the ones in $\mathcal{B}%
_{n}$. Given a graph $G$ on $m\leq n$ vertices which are labeled by disjoint
subsets of $[n]$, and two vertices labeled $s$ and $t$ we obtain $G/\{s,t\}$%
, the \emph{contraction} of $G$ on $\{s,t\}$, where 
\begin{equation*}
V(G/\{s,t\})=V(G)\setminus \{s,t\}\bigcup \{st\}\quad \mbox{and}\quad
E(G/\{s,t\})=E(G)\setminus E(s,t)
\end{equation*}%
where $E(s,t)$ is the set of edges between the vertices $s$ and $t$. We
label the vertices $K_{n}$ with the singletons $\{1\},\ldots ,\{n\}$, and we
call a graph $G$ obtained by a sequence of contractions from $K_{n}$ a
contraction of $K_{n}$. Contractions of $K_{n}$ form a lattice where $H\geq
G $ if H can be obtained by a sequence of contractions from $G$. This
lattice is isomorphic to the partition lattice $\Pi _{n}$ which is in turn
isomorphic to the lattice of flats of $K_{n}$ ordered by inclusion: a flat
of $K_{n}$ is the set of edges that are \emph{not} present in a contraction
of $K_{n}$. Figure \ref{Figure1} illustrates the lattice of contractions of $%
K_{4}$.

\noindent Now let $\mathcal{F} = \{\emptyset = F_0 \subset F_1 \subset F_2
\subset \cdots \subset F_{n-2} \subset F_{n-1} = {\binom{[n] }{2}} \}$ be a
maximal chain of flats of $K_n$ obtained from $K_n$ by a sequence of $n-1$
contractions to $K_1$ with the vertex label $[n]$. Note that $F_i \setminus
F_{i-1}$ is $E(s,t)$ for the corresponding contraction. We define a cone
that is associated to $\mathcal{F}$ as 
\begin{equation*}
\begin{array}{lcl}
C_\mathcal{F} & = & \Big\{ (x(i,j)) \in \mathbb{R}^{\binom{n }{2}} \, : \,
\\ 
&  & \{x(k,l) = x(s,t) : (k,l), (s,t) \in F_1 \setminus F_0\} \quad \leq \\ 
&  & \{x(k,l) = x(s,t) : (k,l), (s,t) \in F_2 \setminus F_1\} \quad \leq
\cdots \leq \\ 
&  & \{x(k,l) = x(s,t) : (k,l), (s,t) \in F_{n-1} \setminus F_{n-2}\} \Big\}.%
\end{array}%
\end{equation*}
The set of $C_\mathcal{F}$ as $\mathcal{F}$ ranges over all maximal chains
in $\Pi_n$ is the maximal cones of $\mathcal{F}_n$. As we mentioned above
the maximal cones of the Bergman complex $\mathcal{B}_n$ are refined by the
cones in $\mathcal{F}_n$. Indeed, two cones $C_{\mathcal{F}^1}$ and $C_{%
\mathcal{F}^2}$ belong to the same maximal cone in $\mathcal{B}_n$ if the
chain of flats $\mathcal{F}^1$ and $\mathcal{F}^2$ differ exactly in one
flat, say $F_i^1 \neq F_i^2$, and $(F_i^1 \setminus F_{i-1}^1) \cap (F_i^2
\setminus F_{i-1}^2) = \emptyset$. \newline

\begin{figure}[tbp]
\centerline{
\includegraphics[height=8cm]{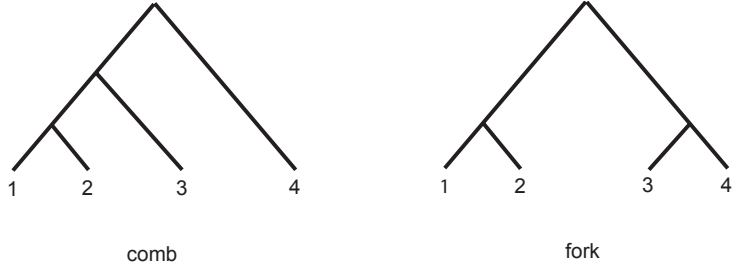}
}
\caption{Comb and fork trees on four leaves}
\label{Figure2}
\end{figure}

\begin{ex}
There are two types of cones corresponding to two types of flag of flats in $%
K_{4}$, namely, 
\begin{equation*}
\mathcal{F}=\{\emptyset \subset \{(1,2)\}\subset
\{(1,2),(1,3),(2,3)\}\subset {\binom{[4]}{2}}\}\mbox{   and}
\end{equation*}%
\begin{equation*}
\mathcal{F}^{\prime }=\{\emptyset \subset \{(1,2)\}\subset
\{(1,2),(3,4)\}\subset {\binom{[4]}{2}}\}.
\end{equation*}%
These go with two types of trees on four leaves: the comb and the fork in
Figure \ref{Figure2}. The corresponding maximal cones in $\mathcal{F}_{4}$
are: 
\begin{equation*}
C_{\mathcal{F}}=\Big\{(x(i,j))\in \mathbb{R}^{6}\,:\,x(1,2)\leq
x(1,3)=x(2,3)\leq x(1,4)=x(2,4)=x(3,4)\Big\}
\end{equation*}%
\begin{equation*}
C_{\mathcal{F}^{\prime }}=\Big\{(x(i,j))\in \mathbb{R}^{6}\,:\,x(1,2)\leq
x(3,4)\leq x(1,3)=x(1,4)=x(2,3)=x(2,4)\Big\}
\end{equation*}
\end{ex}

\begin{prop}
The UPGMA algorithm produces an equidistant tree.
\end{prop}

\begin{proof} It is clear that this algorithm performs a sequence of contractions
 starting from $K_n$ and ending in $K_1$. At iteration $i$ of the repeat loop
 we let $F_i \setminus F_{i-1}$ to be $E(s,t)$ that has been identified. 
 The algorithm sets $x(i,j) = x(k,l)$ for all $(i,j) , (k,l) \in E(s,t)$. 
 So we just need to show that $x(a,b) \leq x(c,d)$ for $(a,b) \in F_i \setminus F_{i-1}$
 and $(c,d) \in F_{i+1} \setminus F_i$.  We let the edges identified in the $(i+1)$st loop
 to be $E(v,w)$. There are two cases: either one of $v$ or $w$ is $s \cup t$  or not.
 In the second case 
 $$x(a,b) \, =\,  \frac{1}{C(s,t)} \sum_{(i,j) \in E(s,t)} d(i,j)  \quad \leq \quad \frac{1}{C(v,w)} \sum_{(i,j) \in E(v,w)} d(i,j) \, = \, x(c,d).$$ 
 For the first case, without loss of generality we assume $v = s \cup t$ and
 $$ \frac{1}{C(s,t)} \sum_{(i,j) \in E(s,t)} d(i,j)  \, \leq \, 
\frac{1}{C(s,w)} \sum_{(i,j) \in E(s,w)} d(i,j) \, \leq \,
\frac{1}{C(t,w)} \sum_{(i,j) \in E(t,w)} d(i,j).$$ 
Now since $\frac{1}{C(v,w)} \sum_{(i,j) \in E(v,w)} d(i,j) $ is equal to 
$$\frac{C(s,w)}{C(v,w)} \left(\frac{1}{C(s,w)} \sum_{(i,j) \in E(s,w)} d(i,j)\right) + \frac{C(t,w)}{C(v,w)}\left(\frac{1}{C(t,w)} \sum_{(i,j) \in E(t,w)} d(i,j)\right)$$
and $C(v,w) = C(s,w) + C(t,w)$
we get the desired inequality in this case as well.
 \end{proof}In the rest of the paper, we will denote the cone in $\mathcal{F}%
_{n}$ which the UPGMA identifies as $C_{UPGMA}$.

\begin{prop}
\label{prop:upgmaproj} If $(x(i,j))$ is the vector that the UPGMA outputs on
the input of the vector $(d(i,j))$ then $(x(i,j))$ is the orthogonal
projection of $(d(i,j))$ onto $C_{UPGMA}$.
\end{prop}

\begin{proof} Let  $\cF = \{\emptyset = F_0 \subset F_1 \subset F_2 \subset \cdots \subset F_{n-2} 
\subset  F_{n-1} = [n]\}$ be the chain of flats  that define the cone $C_{UPGMA}$.
Let $L_{UPGMA}$ be the  smallest subspace containing $C_{UPGMA}$. This subspace is defined by 
 $$ \begin{array}{lcl}
 L_ {UPGMA} & = & \Big\{ (x(i,j)) \in \rr^{n \choose 2} \, : \, \\
 & &  x(k,l) = x(s,t) \,\, \forall \,\,  (k,l), (s,t) \in F_1 \setminus F_0  \quad \mbox{and} \\
 & &  x(k,l) = x(s,t) \,\, \forall \,\,  (k,l), (s,t) \in F_2 \setminus F_1 \quad \mbox{and} \quad \cdots \quad \\
 & &  x(k,l) = x(s,t) \,\, \forall \,\, (k,l), (s,t) \in F_{n-1} \setminus F_{n-2} \Big\},\end{array}$$  
and it has an orthonormal basis consisting of the set of  vectors 
$$\Big\{  \frac{1}{\sqrt{|F_i \setminus F_{i-1}|}} \sum_{(s,t) \in F_i \setminus F_{i-1}} e(s,t)  \quad  : \quad
i = 1, \ldots, n-1 \Big\}
$$
where $e(s,t) \in \rr^{n \choose 2}$ is the standard unit vector corresponding to the edge $(s,t)$ of 
$K_n$.  The linear projection formula with respect to this orthonormal basis 
implies that $x(v,w)$ coordinate of the projection of $(d(i,j))$ is equal to 
$$\frac{1}{|F_k \setminus F_{k-1}|} \sum_{(i,j) \in F_k \setminus F_{k-1}} d(i,j)$$
with $(v,w) \in F_k \setminus F_{k-1}$.  Note that if $(v,w)$ belongs to the contracted 
edges $E(s,t)$ during the UPGMA that produced $\cF$, then
$E(s,t) = F_k \setminus F_{k-1}$ and $C(s,t) = |F_k \setminus F_{k-1}|$, and therefore the 
projected vector $(x(i,j))$ is precisely the vector generated by UPGMA. Therefore this projected
vector is not only in $L_{UPGMA}$ but also in $C_{UPGMA}$. This shows that 
UPGMA performs an orthogonal projection of $(d(i,j))$ onto $C_{UPGMA}$.
 \end{proof}

\begin{cor}
When $n=3$ UPGMA produces the least squares tree.
\end{cor}

\begin{proof} We assume that $d(1,2) \leq d(1,3) \leq d(2,3)$. 
 The two fans $\cB_3$ and $\cF_3$ are identical with three cones 
 described by the three chains of flats
 $$ \cF_{12} = \{ \emptyset \subset \{(1,2)\} \subset {[3] \choose 2}\},  \,\,
 \cF_{13} = \{ \emptyset \subset \{(1,3)\} \subset {[3] \choose 2} \},  \,\,
 \cF_{23} = \{ \emptyset \subset \{(2,3)\} \subset {[3] \choose 2}\}.  $$
 UPGMA produces 
 the tree in $C_{\cF_{12}}$ where the leaves labeled with $1$ and $2$ form a cherry, and
 $$x(1,2) = d(1,2) \quad \mbox{and} \quad x(1,3) = x(2,3) = (d(1,3) + d(2,3))/2.$$
 The square distance of this tree to the data point is $(d(2,3) - d(1,3))^2/2$. 
 We can orthogonally project the data point onto  $L_{\cF_{13}}$
 and $L_{\cF_{23}}$ to obtain 
 $x(1,3) = d(1,3), \,\, x(1,2)=x(2,3) = (d(1,2) + d(2,3))/2$ and 
 $x(2,3) = d(2,3), \,\, x(1,2)=x(1,3) = (d(1,2) + d(1,3))/2$,
 respectively. The first projection is in $C_{\cF_{13}}$ if and only if $d(1,3) \leq (d(1,2)+d(2,3))/2$,
 and the second projection is never in $C_{\cF_{23}}$ unless $d(1,2)=d(1,3)=d(2,3)$. 
 Theorem \ref{thm:main}  implies that the best tree is either the UPGMA tree or the tree obtained
 from $C_{\cF_{13}}$ if the projection falls into this cone.  
 Since the square distance from the data point to this projection  is 
 $(d(2,3) - d(1,2))^2/2 \geq (d(2,3) - d(1,3))^2/2$ we get the result.
 \end{proof}

\begin{ex}
When $n=4$ UPGMA tree may be arbitrarily worse than the least squares tree.
Let the data be $(d(i,j))=(d(1,2),\ldots ,d(3,4))=(1,2,20,10,28+\epsilon ,5)$%
. The UPGMA tree is obtained by contracting the edge $(1,2)$ and then $(3,4)$
in $K_{4}$. This gives us $(x(i,j))=(x(1,2),\ldots ,x(3,4))=(1,15+\frac{1}{4}%
\epsilon ,15+\frac{1}{4}\epsilon ,15+\frac{1}{4}\epsilon ,15+\frac{1}{4}%
\epsilon ,5)$, and the square distance from $(d(i,j))$ to $(x(i,j))$ is $%
388+31\epsilon +\frac{3}{4}\epsilon ^{2}$. The data point can also be
orthogonally projected onto the cone $C_{\mathcal{F}}$ where $\mathcal{F}%
=\{\emptyset \subset F_{1}\subset F_{2}\subset F_{3}={\binom{[4]}{2}}\}$
with $F_{1}=\{(1,2)\}$, $F_{2}\setminus F_{1}=\{(1,3),(2,3)\}$, and $%
F_{3}\setminus F_{2}=\{(1,4),(2,4),(3,4)\}$. The resulting point is $%
(y(i,j))=(1,6,6,\frac{53}{3}+\frac{1}{3}\epsilon ,\frac{53}{3}+\frac{1}{3}%
\epsilon ,\frac{53}{3}+\frac{1}{3}\epsilon )$, and the square distance from $%
(d(i,j))$ to $(y(i,j))$ is $\frac{914}{3}+\frac{214}{9}\epsilon +\frac{2}{3}%
\epsilon ^{2}$. The first expression is greater than the second one for any $%
\epsilon \geq 0$. Indeed, the difference is $\frac{250}{3}+\frac{65}{9}%
\epsilon +\frac{1}{12}\epsilon ^{2}$, and this shows that the UPGMA tree
could be arbitrarily bad.
\end{ex}


\section{The geometry of projections}

\label{sec:projections}

\noindent In the preceding section we showed that UPGMA performs an
orthogonal projection of $(d(i,j))$ onto a distinguished cone of the complex 
$\mathcal{F}_{n}$. It is not immediately clear whether the least squares
equidistant tree is obtained by projecting $(d(i,j))$ \emph{orthogonally}
onto \emph{some} cone of $\mathcal{F}_{n}$. Such a tree will be obtained by
locating a point on $\mathcal{F}_{n}$ (a polyhedral complex) that is closest
to $(d(i,j))$, and in general, nearest point maps of polyhedral complexes do
not have to be given by orthogonal projections onto the \emph{maximal}
faces: take for instance the polyhedral complex in $\mathbb{R}^{2}$ whose
maximal faces are the nonnegative $x$-axis together with the nonnegative $y$%
-axis. For any point with negative coordinates the nearest point is the
origin. Although this is obtained by an orthogonal projection onto the
origin, these projections are not orthogonal to the maximal faces. In this
section, we first show that for $\mathcal{F}_{n}$ and the Bergman complex $%
\mathcal{B}_{n}$ the unexpected happens.

\noindent We start with a definition. Given a maximal chain of flats $%
\mathcal{F}$ of $K_n$ as in Section \ref{sec:bergman} we let $P_\mathcal{F}$
to be the set of points in $\mathbb{R}^{\binom{n }{2}}$ that orthogonally
projects to some point in $C_\mathcal{F}$. Since $P_\mathcal{F} = C_\mathcal{%
F} + L_\mathcal{F}^\perp$ where $L_\mathcal{F}$ is the smallest subspace
containing $C_\mathcal{F}$, it is clear that $P_\mathcal{F}$ is also a
polyhedral cone. We call this cone the \emph{projection cone} of $C_\mathcal{%
F}$.

\begin{thm}
\label{thm:covering} The projection cone $P_\mathcal{F}$ is the
full-dimensional cone defined by the $n-2$ inequalities 
\begin{equation*}
\frac{1}{|F_{k} \setminus F_{k-1}|} \sum_{(i,j) \in F_k \setminus F_{k-1}}
x(i,j) \,\, \leq \,\, \frac{1}{|F_{k+1} \setminus F_{k}|} \sum_{(i,j) \in
F_{k+1} \setminus F_{k}} x(i,j)
\end{equation*}
where $k = 1, \ldots, n-2$. The common refinement of $P_\mathcal{F}$ over
all $\mathcal{F}$ is a complete polyhedral fan.
\end{thm}

\begin{proof} Let $K_\cF$ be the cone defined by the above inequalities. The proof of Proposition \ref{prop:upgmaproj} implies that any point in $K_\cF$
projects to a point in $C_\cF$: one should only note that if $(x(i,j))$ satisfies the
inequalities then $(x(i,j)) + (p(i,j))$ also satisfies them for any $(p(i,j))$ in 
$L_\cF^\perp$ since vectors in $L_\cF^\perp$ do not change the averages which are
on both sides of these inequalities. Conversely, any point in $P_\cF$ is of the form
$(y(i,j)) + (p(i,j))$ where $(y(i,j)) \in C_\cF$ which trivially satisfies these inequalities.
The intersection of any collection of $P_\cF$ is a nonempty cone since 
the intersection of {\em all} $P_\cF$ contains the line generated by $(1,1, \ldots, 1)$ (in fact,
it  is equal to this line). Moreover, Proposition \ref{prop:upgmaproj} implies that
every point in $\rr^{n \choose 2}$ is in {\em some} $P_{UPGMA}$. This shows that the common
refinement of $P_\cF$ is a complete polyhedral fan. 
\end{proof}

\noindent At the end of this section we will look more carefully at this
polyhedral complex obtained by superimposing all $P_{\mathcal{F}}$. For our
main result we need two technical lemmas.

\begin{lemma}
\label{lem:technical1} Suppose $\mathcal{F}^1$ and $\mathcal{F}^2$ are two
distinct maximal chains of flats of $K_n$. Then the interior of $P_{\mathcal{%
F}^1}$ and the cone $C_{\mathcal{F}^2}$ are disjoint.
\end{lemma}

\begin{proof} Suppose  $\cF^j = \{ \emptyset = F_0^j \subset F_1^j \subset \cdots \subset F_{n-2}^j \subset 
F_{n-1}^j = { [n] \choose 2} \}$ for $j =1,2$.  Note that the
relative interior of $P_{\cF^1}$ is defined by the inequalities defining $P_{\cF^1}$ except that
$\leq$ are replaced by $<$. We suppose that the intersection of the
interior of $P_{\cF^1}$ and the cone $C_{\cF^2}$ is not empty, and we will reach a contradiction.
Assume that $F_p^1 = F_p^2$ for $p < q$ and $F_q^1 \neq F_q^2$. Let $(x(i,j))$ be a point in
this nonempty intersection, and let $(y(i,j))$ be the projection of this point onto $C_{\cF^1}$. 
Let $(s,t)$ be an edge in $F_q^1 \setminus F_{q-1}^1$. Then we know that $y(i,j) = a$ for all
$(i,j) \in F_q^1 \setminus F_{q-1}^1$ and therefore $y(s,t) = a$. The edge $(s,t)$ is in 
$F_r^2 \setminus F_{r-1}^2$ where $r > q$ since otherwise $F_q^1 = F_q^2$. Since
$F_r^2$ is a flat containing $F_{q-1}^2=F_{q-1}^1$ the general theory of matroids implies 
that $F_r^2 \setminus F_{r-1}^2$ contains $F_q^1 \setminus F_{q-1}^1$. 
Because $(x(i,j))$ is in $C_{\cF^2}$ we conclude that $x(i,j) =b$ for all
$(i,j) \in F_r^2 \setminus F_{r-1}^2$ and therefore $x(i,j) = b$ for
all $(i,j) \in F_q^1 \setminus F_{q-1}^1$. The orthogonal projection onto $C_{\cF^1}$
keeps the average of these $x(i,j)$ constant. In other words, $a=b$ and $x(i,j) = y(i,j)$
for  $ (i,j) \in F_q^1 \setminus F_{q-1}^1$. Now let $(u,v)$ be an edge in $F_{q+1}^1 \setminus F_q^1$.
Again we know that $y(i,j) = c > a$ for all $(i,j) \in F_{q+1}^1 \setminus F_q^1$, including 
$y(u,v) = c$. 
We will show that $(u,v)$  is in $F_z^2 \setminus F_{z-1}^2$ where $z > r$. Suppose not.
Then $F_r^2 \setminus F_{q-1}^2$ contains $F_{q+1}^1 \setminus F_{q-1}^1$, and
this implies that $x(i,j) = b_{i,j} \leq a$ for all $(i,j) \in F_{q+1}^1 \setminus F_{q-1}^1$. But then
the orthogonal projection argument implies that $y(i,j) = c \leq a$ for $(i,j) \in F_{q+1}^1 \setminus F_q^1$.
This is a contradiction and we conclude that $z>r$. The above chain of arguments can be applied
to all $F_{k}^1 \setminus F_{k-1}^1$ for $k = q, \ldots, n-1$ to produce a chain  $
F_{r_q}^2 \subset F_{r_{q+1}}^2 \subset \cdots \subset F_{r_{n-1}}^2$
where $  q < r_q < r_{q+1} < \cdots < r_{n-1} \leq n-1$ (we have constructed the first two
members of this chain, namely $r_q = r$ and $r_{q+1} = z$). However, this is
a contradiction since there are only $n-1-q$ distinct integers bigger than $q$ and at most $n-1$.
\end{proof}  

\begin{lemma}
\label{lem:technical2} Suppose $\mathcal{F}^1$ and $\mathcal{F}^2$ are two
distinct maximal chains of flats in $K_n$. If $P_{\mathcal{F}^1} \cap C_{%
\mathcal{F}^2}$ is nonempty, then this intersection is contained in $C_{%
\mathcal{F}^1} \cap C_{\mathcal{F}^2}$.
\end{lemma}

\begin{proof} This proof invokes similar ideas as in the proof of Lemma \ref{lem:technical1}.
Suppose  $\cF^j = \{ \emptyset = F_0^j \subset F_1^j \subset \cdots \subset F_{n-2}^j \subset 
F_{n-1}^j = { [n] \choose 2} \}$ for $j =1,2$. Assume that $F_p^1 = F_p^2$ for $p < q$ and $F_q^1 \neq F_q^2$. Let $(x(i,j))$ be a point in $P_{\cF^1} \cap C_{\cF^2}$, and let $(y(i,j))$ be the projection of 
this point onto $C_{\cF^1}$. We will show that $x(i,j) = y(i,j)$ for all $(i,j)$. By our assumption 
this is true for  all $(i,j) \in F_{q-1}^1 = F_{q-1}^2$. 
Let $(s,t)$ be an edge in $F_q^1 \setminus F_{q-1}^1$. Then we know that $y(i,j) = a$ for all
$(i,j) \in F_q^1 \setminus F_{q-1}^1$ and therefore $y(s,t) = a$. The edge $(s,t)$ is in 
$F_r^2 \setminus F_{r-1}^2$ where $r > q$ since otherwise $F_q^1 = F_q^2$. Since
$F_r^2$ is a flat containing $F_{q-1}^2=F_{q-1}^1$ we conclude
that $F_r^2 \setminus F_{r-1}^2$ contains $F_q^1 \setminus F_{q-1}^1$. 
Because $(x(i,j))$ is in $C_{\cF^2}$ we have $x(i,j) =b$ for all
$(i,j) \in F_r^2 \setminus F_{r-1}^2$ and therefore $x(i,j) = b$ for
all $(i,j) \in F_q^1 \setminus F_{q-1}^1$. The orthogonal projection onto $C_{\cF^1}$
keeps the average of these $x(i,j)$ constant. In other words, $a=b$ and $x(i,j) = y(i,j)$
for  $ (i,j) \in F_q^1 \setminus F_{q-1}^1$. Now let $(u,v)$ be an edge in $F_{q+1}^1 \setminus F_q^1$.
Now we know that $y(i,j) = c \geq a$ for all $(i,j) \in F_{q+1}^1 \setminus F_q^1$, including 
$y(u,v) = c$. 
Furthermore  $(u,v)$ is in $F_z^2 \setminus F_{z-1}^1$ where $z > q-1$. If $z \leq r$ then 
$F_{q+1}^1 \setminus F_{q-1}^1$ is contained in $F_r^2$, and therefore $x(i,j) = b_{i,j} \leq a$ for
all $(i,j) \in F_{q+1}^1 \setminus F_q^1$. If $x(i,j)<a$ for any of these edges, then
the average of $x(i,j)$ in this set is strictly less than $a$. But this average is equal
to the average of $y(i,j)$ for $(i, j) \in F_{q+1}^1 \setminus F_q^1$, and we get a contradiction
since this implies $c<a$. Therefore, if $z \leq r$ then $a=c$ and $x(i,j) = y(i,j)$ for 
all $(i,j)$ in $F_{q+1}^1 \setminus F_q^1$. If $z >r$, then $F_z^2 \setminus F_{z-1}^2$
contains $F_{q+1}^1 \setminus F_q^1$. Therefore $x(i,j) = b$ for all $(i,j) \in F_{q+1}^1 \setminus F_q^1$,
and the average of these $x(i,j)$ is $b$. But $b = c$ since the average of $x(i,j)$ and $y(i,j)$ is constant
for $(i,j) \in F_{q+1}^1 \setminus F_q^1$. But this implies $x(i,j) = y(i,j)$ for edges in this set
as well. Now we can repeat the same argument for the rest of $F_k^1 \setminus F_{k-1}^1$.

\end{proof}

\begin{thm}
\label{thm:main} Let $\mathcal{P}_d$ be the set of projection cones
containing the data point $(d(i,j))$. Then the best least squares
equidistant tree corresponds to a point in $C_\mathcal{F}$ for some $P_%
\mathcal{F} \in \mathcal{P}_d$.
\end{thm}

\begin{proof}
Let $x = (x(i,j))$ be the point corresponding to the best least square equidistant tree, and 
let $C$ be the cone of $\cF_n$ where $x \in C$. We let $P$ be the projection cone of $C$. 
If the line segment $\overline{x-d}$ is orthogonal to $C$, then $P \in \mathcal{P}_d$ and we
are done. If not, we show that $x$ is on the boundary of $C$ (and hence $P$). 
Suppose that $x$ is in the relative interior of $C$, and hence in the interior of $P$. 
If $\overline{x-d}$ is entirely contained in $P$, then it must be perpendicular to $C$ and this
would imply that $P \in \mathcal{P}_d$. Hence let $y  = (y(i,j))$
be the first point on $\overline{x-d}$ which intersects $P$. Clearly, $y$ must
be on the boundary of $P$. If we let $y'$ be the orthogonal projection
of $y$ onto $C$, then we conclude that $|| y' - d || <  || y- d|| + ||y'-y|| <  ||y-d|| + ||x-y|| = ||x-d||$,
and this is a contradiction. Hence $x$ is on the boundary of $C$ and $P$. 
Now let $\mathcal{P}_x$ be the projection cones containing $x$. Since
$x$ is on the boundary of $P \in \mathcal{P}_x$, and by  Lemma \ref{lem:technical1},
none of the projection cones in $\mathcal{P}_x$ can contain $x$ in their interiors. If we
let $\mathcal{P}_x = \{ P_{\cF_1}, \ldots, P_{\cF_k}\}$, then by Lemma \ref{lem:technical2} $x$ is
contained in (the boundaries of) $C_{\cF_1}, \ldots, C_{\cF_k}$. Now only two things can happen.
Either $\overline{x-d}$ is entirely contained in one of $P_{\cF_i}$ where $i=1,\ldots, k$, in which
case this cone also belongs to $\mathcal{P}_d$, and we are done, or otherwise for each $P_{\cF_i}$
there is a first point $y_i$ on $\overline{x-d}$ intersecting $P_{\cF_i}$. Repeating the above argument
we can conclude that $y_i = x$ for all $i$. This means that for every point $y$ on $\overline{x-d}$
$\mathcal{P}_y$ is disjoint from $\mathcal{P}_x$. But since the projection cones are closed cones 
this is a contradiction, unless $x=d$.
\end{proof}

\noindent In the light of Theorem \ref{thm:main} we introduce a graph $%
\mathcal{G}_{n,d}$ associated to each data point $d = (d(i,j))$ in $\mathbb{R%
}^{\binom{n }{2}}$. The vertices of this graph are $C_\mathcal{F}$ where $P_%
\mathcal{F} \in \mathcal{P}_d$, and there is an edge between two vertices $%
C_{\mathcal{F}^1}$ and $C_{\mathcal{F}^2}$ if these two cones in $\mathcal{F}%
_n$ share a facet.

\begin{prop}
The graph $\mathcal{G}_{3,d}$ is either $K_1$, $K_2$ or $K_3$, and when $n
\geq 4$, $\mathcal{G}_{n,d}$ could have more than one component.
\end{prop}

\begin{proof}
At most two out of the three inequalities 
$d_{12} \leq (d_{13} + d_{23})/2$, $d_{13} \leq (d_{12} + d_{23})/2$,  and $d_{23} \leq (d_{12} + d_{13})/2$ hold unless $d_{12} = d_{13} = d_{23}$. In the latter case $\mathcal{G}_d = K_3$, 
and otherwise $\mathcal{G}_d = K_j$ if $j=1,2$ of the inequalities are satisfied.
We use the following data to illustrate that $\mathcal{G}_{4,d}$ can be disconnected:
$$ (d(1,2), d(1,3), d(1,4), d(2,3), d(2,4), d(3,4))  \, = \, (1,2,3,2,7,3).$$
The set $\mathcal{P}_d$ consists of four cones, and hence there are four trees one
can obtain. The component of the UPGMA tree has a total of two trees, and there
are two more components where each component is just one tree. 
\end{proof}

\noindent
We finish this section by studying the polyhedral complex we defined in
Theorem \ref{thm:covering}, namely the common refinement of the projection
cones $P_\mathcal{F}$ for each maximal cell $C_\mathcal{F}$ in the Bergman
complex $\mathcal{B}_n$. We denote this complex by $\mathcal{Q}_n$. Note
that $\mathcal{Q}_n$ is full dimensional complex in $\mathbb{R}^{\binom{n }{2%
}}$, and the interior of a full dimensional cell in $\mathcal{Q}_n$ consists
of those data vectors which project to the same set of $C_\mathcal{F}$.

\begin{ex}
When $n=3$ the complex $\mathcal{Q}_3$ is easy to describe. There are total
of six maximal cells which are of two different types. The first type
consists of those vectors which project to exactly one of the three $C_%
\mathcal{F}$. The second type consists of those vectors which project to
exactly two $C_\mathcal{F}$.
\end{ex}

\begin{ex}
When $n=4$ at most six distinct projection cones could have an intersection
that gives a maximal cell in $\mathcal{Q}_4$ as we checked with a short
MAPLE program. There are a total of $166$ such cells, but they come in ten
different orbits with respect to the action of $S_4$. The following table
lists a representative of each orbit. Since the projection cones are indexed
by binary trees on four leaves, we just list these trees. %
\renewcommand\arraystretch{0.7} \fontsize{10}{12} 
\begin{equation*}
\begin{array}{|c|c|}
\hline
\mbox{orbit representative} & \mbox{orbit size} \\ \hline
(((1,2),3),4) \, (((1,2),4),3) \, (((1,3),2),4) \, (((1,3),4),2) \,
(((1,4),2),3) \, (((1,4),3),2) & 4 \\ 
(((1,2),3),4) \, (((1,2),4),3) \, (((1,3),2),4) \, (((1,3),4),2) \,
(((1,4),2),3) \, ((1,4),(2,3)) & 24 \\ 
(((1,2),3),4) \, (((1,2),4),3) \, (((1,3),2),4) \, (((1,4),2),3) \,
((1,3),(2,4)) \, ((1,4),(2,3)) & 12 \\ 
(((1,2),3),4) \, (((1,2),4),3) \, (((1,3),2),4) \, (((1,4),3),2) \,
(((2,4),1),3) \, ((1,3),(2,4)) & 24 \\ 
(((1,2),3),4) \, (((1,2),4),3) \, (((1,3),2),4) \, (((2,3),4),1) \,
(((2,4),3),1) \, ((1,3),(2,4)) & 24 \\ 
(((1,2),3),4) \, (((1,2),4),3) \, (((1,3),2),4) \, (((2,4),1),3) \,
((1,3),(2,4)) \, ((1,4),(2,3)) & 12 \\ 
(((1,2),3),4) \, (((1,3),2),4) \, (((1,4),2),3) \, (((2,4),1),3) \,
((1,3),(2,4)) \, ((1,4),(2,3)) & 24 \\ 
(((1,2),3),4) \, (((1,3),2),4) \, (((2,4),1),3) \, (((3,4),1),2) \,
((1,2),(3,4)) \, ((1,3),(2,4)) & 12 \\ 
(((1,2),3),4) \, (((1,3),2),4) \, (((2,4),1),3) \, (((3,4),2),1) \,
((1,2),(3,4)) \, ((1,3),(2,4)) & 24 \\ 
(((1,2),3),4) \, (((1,3),2),4) \, (((2,4),3),1) \, (((3,4),2),1) \,
((1,2),(3,4)) \, ((1,3),(2,4)) & 6 \\ \hline
\end{array}%
\end{equation*}
\end{ex}

\begin{thm} There is a maximal cell in $\mathcal{Q}_n$ which is the intersection of at least
$(n-1)!$ projection cones; i.e., there are data vectors in $\mathbb{R}^{{n \choose 2}}$ which
orthogonally project onto at least $(n-1)!$ (non-degenerate) equidistant trees.
\end{thm}
\begin{proof} Let $a < b $ two real numbers and let $x(i,j) \in \mathbb{R}^{{n \choose 2}}$ be the
data vector where $x(1,j) = a$ for $j = 2, \ldots, n$ and $x(i,j) = b$ for all other components. We claim
that this vector is in the interior of the intersection of $(n-1)!$ projection cones corresponding to 
the comb trees of the form $(\cdots ((1,a_2), a_3) \cdots), a_n)$ where $a_2, a_3, \ldots, a_n$ run
through all permutations of $\{2, \ldots, n\}$. For any one of these trees our data vector is in the interior
of the corresponding projection cone if and only if 
\[ a < \frac{a+b}{2} < \frac{a+2b}{3} <  \frac{a+3b}{4} <  \cdots <  \frac{a+(n-2)b}{n-1}.\]
The above inequalities hold for the choice of $a$ and $b$ we made. This proves the theorem.
\end{proof}

\noindent
This theorem asserts that there are maximal cells in $\mathcal{Q}_n$ that are intersections 
of {\em at least} $(n-1)!$ projection cones. We believe that the number of such cones cannot 
exceed $(n-1)!$, though we do not have a proof. 
\begin{conj} 
The maximal cells in $\mathcal{Q}_n$ are obtained as the intersection of at most $(n-1)!$ projection cones.
\end{conj}

\section{Extended UPGMA and  Branch-and-Bound}

\label{sec:comp}

\noindent In view of the results in Section \ref{sec:projections} we propose
two algorithms. The first one is an extension of the usual UPGMA which
searches the component of the graph $\mathcal{G}_{n,d}$ to which the UPGMA
tree belongs to. Even when this component is large this extended UPGMA
algorithm performs well and finds the best tree in this component. The
drawback of this algorithm is that it may not produce the best tree. Our
second algorithm is an exact algorithm which produces the best equidistant
tree with a branch and bound approach on the space of maximal chains of the lattice of 
contractions of $K_n$.  We will present this as a shortest path algorithm
on the Hasse diagram of this lattice. Recall that this
lattice is isomorphic to the partition lattice $\Pi_n$ where maximal chains
 are in bijection with the maximal cones in $\mathcal{F}_n$.

\begin{algorithm}
\label{alg:ext-upgma} Extended UPGMA \newline
\textbf{Input} : Complete graph $K_n$ with edge weights $d(i,j)$. \newline
\textbf{Output}: An equidistant tree $T$ with leaves $X = \{1, \ldots, n\}$
and $x(i,j)$ for each $i,j \in X$. \newline

\vskip 0.3cm \noindent Using Algorithm \ref{alg:upgma} find the UPGMA tree $%
T_{UPGMA}$ and the corresponding cone $\mathcal{C}_{UPGMA}$ in the Bergman
compex $\mathcal{B}_n$. \newline
Let $\mathtt{Visited} := \{T_{UPGMA}\}$, $\mathtt{Active} := \{ T_{UPGMA} \} 
$, and $T_{best} := T_{UPGMA}$. \newline
\textbf{while $\mathtt{Active} \neq \emptyset$ do}

Let $T \in \mathtt{Active} $ and $\mathtt{Active} := \mathtt{Active}
\setminus \{T\}$.

\textbf{for} each $\mathcal{C}_{T^{\prime }} \in \mathcal{B}_n$ which shares
a facet with $\mathcal{C}_T$ \textbf{do}

$\,\,$ \textbf{if} $\mathcal{C}_{T^{\prime }} \in \mathcal{P}_d$ and $%
T^{\prime }\not \in \mathtt{Visited}$ \textbf{then}

$\,\,$ $\mathtt{Active} := \mathtt{Active} \cup \{T^{\prime }\}$ and $%
\mathtt{Visited} := \mathtt{Visited} \cup \{T^{\prime }\}$

$\,\,$ If $\sum (d(i,j) - x_{T^{\prime }}(i,j))^2 \quad < \quad \sum (d(i,j)
- x_{T_{best}}(i,j))^2$ then $T_{best} := T^{\prime }$.

$\,\,$ \textbf{end if}

\textbf{end for} \newline
\textbf{end while} \newline

\vskip 0.2cm \noindent Output $T_{best}$ and $x_{T_{best}}(i,j)$ $1 \leq i <
j \leq n$.
\end{algorithm}

\vskip 1cm \noindent A few remarks about Algorithm \ref{alg:ext-upgma} are
in order: This algorithm searches the component of $\mathcal{G}_{n,d}$ to
which $\mathcal{C}_{UPGMA}$ belongs, and it outputs the best equidistant
tree in this component. If $\mathcal{G}_{n,d}$ consists of a single
component then the algorithm's output is the optimal tree. The search
depends on the following characterization of Ardila and Klivans \cite{AK}
when $\mathcal{C}_{\mathcal{F}^1}$ and $\mathcal{C}_{\mathcal{F}^2}$ share a
facet in $\mathcal{B}_n$. Finally, checking whether a cone $\mathcal{C}_T$ belongs to $%
\mathcal{P}_d$ is trivial by Theorem \ref{thm:covering}.

\begin{prop}
\label{prop:AK} Two maximal cones $\mathcal{C}_{\mathcal{F}^1}$ and $%
\mathcal{C}_{\mathcal{F}^2}$ share a facet in $\mathcal{B}_n$ if and only if
there exists $0 < j < n-1$ such that $F_i^1 = F_i^2$ for all $i=0, \ldots,
n-1$ except $F_j^1 \neq F_j^2$ and $(F_j^1 \setminus F_{j-1}^1) \cap (F_j^2
\setminus F_{j-1}^2) \neq \emptyset$.
\end{prop}

\noindent Our exact algorithm is a modified shortest path algorithm
performed on the Hasse diagram of the partition lattice $\Pi_n$. We first
introduce some notation for this algorithm. We will represent this Hasse
diagram as a directed graph where the nodes are labeled by flats of $K_n$,
and the edges are directed from the minimum element (corresponding to the
empty flat) to the top element (corresponding to the flat $[n] = \{1,\ldots,
n\}$). For each node (flat) $F$ we will keep track of \emph{incoming} edges $%
I_F$ and \emph{outgoing} edges $O_F$. Each edge is directed from a flat $F_i$
to a flat $F_{i+1}$ of next rank such that $F_i \subset F_{i+1}$. Each such
edge $e$ will have two associated numbers, $x(e)$ and $\ell(e)$, which will
be defined throughout the algorithm using the given data $(d(i,j))$.

\begin{algorithm}
\label{alg:exact} Exact least squares \newline
\textbf{Input} : Complete graph $K_n$ with edge weights $d(i,j)$. \newline
\textbf{Output}: The best least square equidistant tree $T$ with leaves $X =
\{1, \ldots, n\}$ and $x(i,j)$ for each $i,j \in X$. \newline

\vskip 0.3cm \noindent Set $\mathtt{Active}_0 := \{ \emptyset \} $, $%
I_\emptyset := \{g\}$, $x(g) := -\infty$, and $\ell(g) := 0$. \newline
Set $V := \mathtt{Active}_0$ and $A := \{\}$. \newline
\textbf{for} $k = 0, 1, \ldots, n-1$ \textbf{do} \newline
$\mathtt{Active}_{k+1} := \{ \}$.

\textbf{while $\mathtt{Active}_k \neq \emptyset$ do}

$\,\,$ Let $F \in \mathtt{Active}_k $ and $\mathtt{Active}_k := \mathtt{%
Active}_k \setminus \{F\}$.

$\,\,$ \textbf{for} each $e=(F,F^{\prime }) \in O_F$ \textbf{do}

$\,\,\,\,\,$ Set $x(e) := \frac{1}{|F^{\prime }\setminus F|} \sum_{(i,j) \in
F^{\prime }\setminus F} d(i,j)$ and $E := \{f \in I_F \, : \, x(e) \geq x(f)
\}$.

$\,\,\,\,\,$ \textbf{if} $E = \emptyset$ \textbf{then} $x(e) := +\infty$ 
\textbf{else} $h := \mathrm{argmin} \{\ell(f) \, : \, f \in E\}$ \textbf{end
if}

$\,\,\,\,\,$ $\ell(e) = \ell(h) + w(e)$ where $w(e) = \sum_{(i,j) \in
F^{\prime }\setminus F} (x(e) - d(i,j))^2$

$\,\,\,\,\,$ \textbf{if} $x(e) < + \infty$ \textbf{then}

$\,\,\,\,\,\,\,\,$ $\mathtt{Active}_{k+1} := \mathtt{Active}_{k+1} \cup
\{F^{\prime }\}$ and $A := A \cup \{e\}$

$\,\,\,\,\,$ \textbf{end if}

$\,\,$ \textbf{end for}

\textbf{end while} \newline
$V := V \cup \mathtt{Active}_{k+1}$. \newline
\textbf{end for}

\vskip 0.1cm \noindent Find the shortest path $P$ from $\emptyset$ to $[n]$
in the graph $G=(V,A)$ with edge weights $w(e)$ for $e \in A$.

\vskip 0.2cm \noindent Output the tree $T$ corresponding to $P$ and $%
x_{T}(i,j)$ $1 \leq i < j \leq n$.
\end{algorithm}

\vskip 0.5cm \noindent \emph{Proof of Correctness}: Each path $P$ in $G$ from the empty flat
to the full flat 
corresponds to a flag $\mathcal{F}$ and hence a cone $\mathcal{C}_\mathcal{F}
$. By the construction of $G$, the $x(e)$ for the edges $e$ on such a path
give a point in $\mathcal{C}_\mathcal{F}$, and this point is the orthogonal
projection of $(d(i,j))$ onto $\mathcal{C}_\mathcal{F}$. In other words, a
path $P$ in $G$ corresponds to a cone $\mathcal{C}_\mathcal{F} \in \mathcal{P%
}_d$. Since $\sum_{e \in P} w(e)$ is the Euclidean distance from the data
point to the projection in $\mathcal{C}_\mathcal{F}$, Theorem \ref{thm:main}
implies the correctness of the algorithm if for each $\mathcal{C}_\mathcal{F}
\in \mathcal{P}_d$ there is a path $P$ in $G$ from $\emptyset$ to $[n]$. We
show by induction on $i$ that $G$ contains the edges $e_i = (F_{i-1}, F_i)$
corresponding to the flag $\mathcal{F}$. It is trivial to check that $e_1 =
(\emptyset, F_1)$ is in $G$. Moreover $x(e_1) = d(i,j)$ where $F_1 =
\{(i,j)\}$. We assume that $e_k$ for $k \leq i-1$ are in $G$. Note that each 
$e_k$ is added to $G$ during the $k$th pass of the outermost for loop. Now $%
x(e_i) = \frac{1}{|F_i \setminus F_{i-1}|} \sum_{(i,j) \in F_i \setminus
F_{i-1}} d(i,j)$, and because $\mathcal{C}_\mathcal{F} \in \mathcal{P}_d$ we
conclude that $x(e_i) \geq x(e_{i-1})$. Since $e_{i-1} \in I_{F_{i-1}}$ the
set $E$ during the pass of the innermost for loop corresponding to $e_i$ is
nonempty and hence $x(e_i)$ stays finite. This means $e_i$ is added to $G$.
\hfill $\Box$

\vskip 0.3cm For the purposes of the exposition of Algorithm \ref{alg:exact}
we have chosen to first construct the graph $G$ in the algorithm and then
solve the shortest path problem on this graph. In fact one can skip the
construction of $G$ if one adds a pointer to each edge $e$ that points to
the corresponding edge $h$ in the algorithm. With these pointers one can
reconstruct the shortest path and hence the best equidistant tree $T$ at the
end of the algorithm. Note that this algorithm is a branch and bound algorithm
on the space of  all maximal chains in $\Pi_n$ starting from the empty flat: whenever $x(e) = +\infty$ for some $e$ being considered then all such maximal chains containing $e$ are pruned from the branch and bound tree. The branching step is realized when we extend a chain terminating at the node labeled
$F$ by adding $e = (F,F') \in O_F$ for all such edges where $x(e) < +\infty$.


\section{A biology example}

\label{sec:biology} \noindent How does the least squares approach compare to
Bayesian and maximum likelihood methods in practise? \ We compared the\
different methods on a problem in evolution for which some of the data shows
clock-like behavior. \ Murphy et al. (2001) have studied the timing
and sequence of appearence of the mammalian orders using a large DNA
database that includes 42 placental mammals from all orders, plus two
marsupials as the outgroup. \ The model of sequence evolution employed by
Murphy et al. (and by us) was the general-time-reversible+$\Gamma $%
+invariants model. \ Bayesian and maximum likelihood methods converged on
the same combinatorial type of tree (Murphy et al., 2001). \ Distances
estimated during likelihood fitting using this model do not satisfy the
clock hypothesis over the complete dataset; however a subset of eleven
species do show clock-like substitution rates (Murphy et al.,2001,
supplemental material). \ Distances from ten of these taxa were analyzed
here using the exact least squares algorithm. \ The main conclusions of
Murphy et al. on the branching sequence are supported by the best least
squares tree: \ first the Afrotherians, then the Xenartharns and finally the
Boreoeutherians separate from their placental ancestors. \ The only
difference between the ten taxa least squares and maximum likelihood trees
is the position of the dolphin. \ Murphy et al. scaled their tree to obtain
dates using 50 mya for the cat/canid divergence. \ Scaling the best least
squares tree in the same way gives 107 and 101 million years ago for the
bifurcations producing the Afrotherians and Xenarthrans, respectively. \ The
corresponding values reported by Murphy et al., 2001, are 103 and 95 million
years ago. \ Hence the agreement between the least squares and Bayesian or
likelihood methods is quite good.

\begin{figure}[tbp]
\centerline{
\includegraphics[height=12cm]{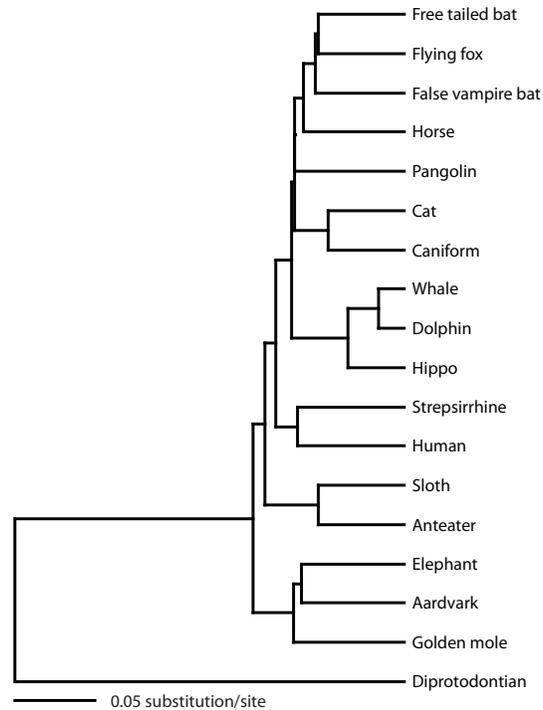}
}
\caption{Example of mammalian phylogeny obtained by extended UPGMA}
\label{Figure3}
\end{figure}

Visual inspection of the complete phylogram from the 44 taxa dataset
suggested that others were nearly contemporaneous with those in the eleven
taxa subset. \ For a sequence of datasets ranging from eleven to nineteen
species, three trees were identified: \ we found the maximum likelihood
tree, the maximum likelihood equidistant tree, and the best least squares
equidistant tree. \ Species added to the eleven taxa subset were the
roussette fruit bat, anteater, whale, hippopotamus, aardvark, human, horse
and sciurid. \ The inexact form of the least squares equidistant tree
algorithm was used on these datasets with more than ten taxa. \ For the
trees with twelve up to eighteen taxa the number of possible least squares
trees was either one or two, with the best least squares tree being the
UPGMA tree in each case. \ The eighteen taxa dataset had two possible trees,
the better was the non-UPGMA tree (Figure 3). \ When the sciurid data was
then added to create a dataset with nineteen taxa, the number of possible
trees jumped to six. 

The number of possible trees for datasets up
to eighteen taxa is small compared to the conjectured $(n-1)!$ upper limit
of trees, indicating that the distances were close to clock-like. \ Hence
the corresponding equidistant trees should be good approximations to the
phylogeny. \ As distances that deviate more from clock-like behavior are
added,  the number of possible trees increases and the
equidistant tree gives a poorer account of the phylogeny.

When two least squares equidistant trees were possible for a given dataset,
the oldest bifurcations were conserved between the two, with the differences
appearing in more recent branchings. \ This observation is expected. \ For a
given internal node the distance to a leaf is one-half the average of all
path lengths between pairs of leaves that pass through that node. \ More
paths pass through the older nodes, so their ages are estimated more
accurately. \ Unless old bifurcations occur very close to each other, they
will be more stable in the set of possible trees. \ The best least squares
trees with up to eighteen taxa all confirmed the branching order of the
Afrotherians, Xenarthrans and Boreoeutherians observed by Murphy et al.,
2001.

There was one persistent difference between the likelihood and least squares
approaches: \ the equidistant least squares trees placed the cetartiodactyls
as an outgroup to the carnivores, bats and pangolin, whereas the maximum
likelihood trees put the bats as an outgroup. \ It is a bit surprising,
since the likelihood and distance methods are both consistent in the
statistical sense, and therefore expected to converge on the same, correct,
tree (Felsenstein, 2004). \ The dataset contains 17028 characters, but
perhaps more data is needed, or a different sample of sequences, to obtain
convergence on one tree.

\vskip 0.5cm 
\noindent
{\bf\large Acknowledgements} The research presented in this work was supported by 
NSF-UBM EF-0436313.

\vfill 
\eject



\end{document}